\documentclass{article}

\usepackage{pifont}
\usepackage{graphics}
\usepackage{enumerate}

\newcommand{\proc}[1]{\noindent\textit{#1}}

\newcommand{\abs}[1]{\left| #1 \right|}

\newcommand{\quotient}[2]{#1/#2}
\newcommand{\set}[2]{\left\{ #1 \, \left| \, #2 \right.\right\}}

\renewcommand{\hat}{\widehat}

\newcommand{\deff}[3]{#1\colon #2 \rightarrow #3}

\usepackage{bbm}
\usepackage{dsfont}
\usepackage{bbold}
\usepackage{amsmath,amssymb}
\usepackage{tabularx}
\usepackage{geometry}
\usepackage{microtype}
\usepackage{tikz}
\usetikzlibrary{calc,arrows,decorations.pathmorphing,backgrounds,positioning,fit,shapes.arrows,shapes.geometric,decorations.pathreplacing} 
\usepackage{caption}
\usepackage{subcaption}
\usepackage{subfig}
\usepackage{bbm}
\usepackage{dsfont}            
\usepackage{float}
\usepackage{bbold}
\RequirePackage{graphics}
\usepackage{amsmath,amsfonts,amssymb,amsthm}

\theoremstyle{plain}

\newtheorem{theorem}{Theorem}[section]
\newtheorem{proposition}[theorem]{Proposition}

\newtheorem{corollary}[theorem]{Corollary}

\newtheorem*{theorem*}{THEOREM}
\newtheorem*{lemma*}{Sublemma}

\theoremstyle{definition}
\newtheorem{definition}[theorem]{Definition}

\theoremstyle{remark}
\newtheorem{remark}[theorem]{Remark}
\newtheorem*{remark*}{Remark}

\begin{document}
\title{Finite-rank Bratteli-Vershik diagrams are expansive -- a new proof}
\author{Siri-Mal\'en H\o ynes, Norwegian University of Science and Technology (NTNU)}
\maketitle
\begin{abstract}
In \cite{DM}, Downarowicz and Maass proved that the Cantor minimal system associated to a properly ordered Bratteli diagram of finite rank is either an odometer system or an expansive system. 
We give a new proof of this truly remarkable result which we think is more transparent and easier to understand.  We also address the question (QUESTION 1) raised in \cite{DM} and we find a better
(i.e.\ lower) bound than the one given in \cite{DM}. In fact, we conjecture that the bound we have found is optimal.
\end{abstract}

\section{Introduction.}
The aim of this paper is to give a new proof of the following result.

\begin{theorem*}
 Let $(V,E,\geq)$ be a properly ordered Bratteli diagram, and let $(V,E)$ be of finite rank. Then the associated Bratteli-Vershik system $(X_{(V,E)}, T_{(V,E)})$ is either an odometer system or an expansive
system.
\end{theorem*}

\begin{remark}
 It is well known that an expansive Cantor minimal system is (conjugate to) a minimal subshift on a finite alphabet (cf.\ Proposition \ref{prop:subshift}). The THEOREM implies that if $(V,E\geq)$ is a properly ordered Bratteli diagram of finite rank \emph{and} $(V,E)$ has the ERS-property (cf.\ Section \ref{sec:bratteli}), then $(X_{(V,E)},T_{(V,E)})$ is either an odometer or a Toeplitz flow \cite{hoynes}.
\end{remark}

In our judgment the proof given of the THEOREM in \cite{DM} is not easy to follow, so we feel that a more transparent proof -- thus hopefully making it more accessible -- is in order for such an important and, frankly speaking, rather surprising
result. We also address the question (QUESTION 1) that is raised in \cite{DM} about finding a better (i.e.\ lower) bound than the one they give in their ``Infection Lemma'', and we do indeed find a significantly lower bound
(cf.\ Corollary \ref{cor:bound}), which we conjecture is optimal. 

We will adopt some of the definitions and terminology from \cite{DM}, but in contrast to \cite{DM} we interpret the definitions directly in terms of the properly ordered Bratteli diagram in question. We 
feel this makes it much easier to grasp the contents of the various definitions. There is also too much ``hand-waving'', we feel, in \cite{DM}, which results in some questionable short-cuts and claims.
But we hasten to say that these can all be fixed by a more careful construction. (Cf.\ Remark at the end of this paper, where we are more specific on this point.) However, we strongly emphasize
that our proof is very much inspired and motivated by the proof in \cite{DM}.

\section{Bratteli diagrams and Bratteli-Vershik systems.}
\label{sec:bratteli}
\begin{figure}
\centering
 \includegraphics{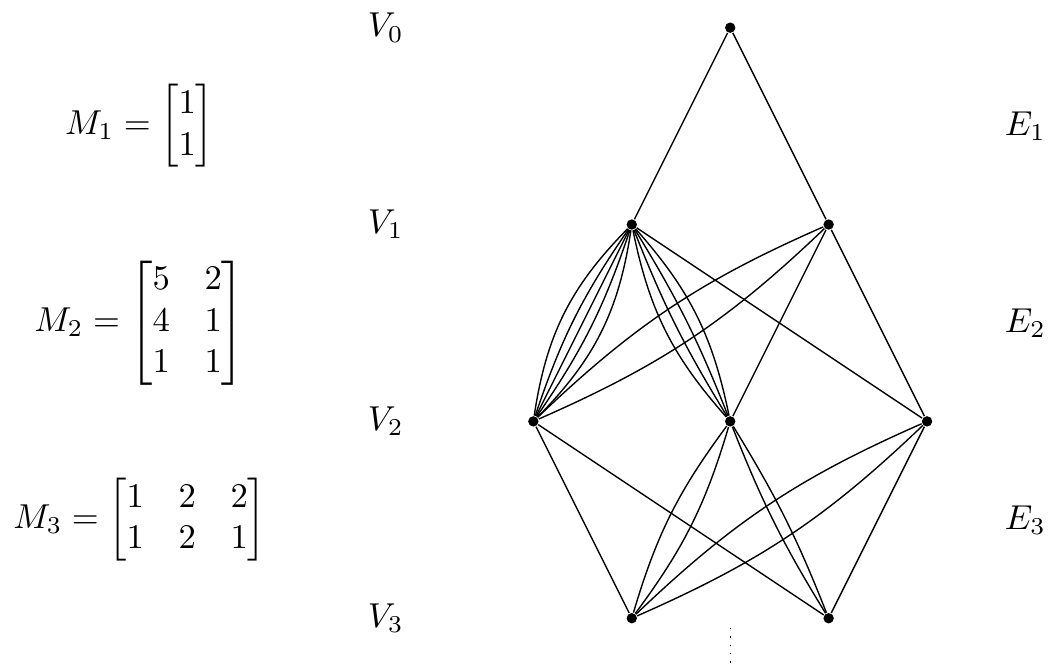}
\caption{An example of a Bratteli diagram}
\label{fig:bD}
\end{figure}

General references for this section are \cite{HPS} and \cite[Section 3]{GPS}.
A Bratteli diagram $(V,E)$ consists of a set of vertices $V= \sqcup_{n=0}^\infty V_n$ and a set of edges $E=\sqcup_{n=1}^\infty E_n$, where the $V_n$'s and the $E_n$'s are finite disjoint sets and where
$V_0=\{v_0\}$ is a one-point set. The edges in $E_n$ connect vertices in $V_{n-1}$ with vertices in $V_n$. If $e$ connects $v \in V_{n-1}$ with $u \in V_n$ we write $s(e)=v$ and $r(e)=u$, where $\deff{s}{E_n}{V_{n-1}}$ and
$\deff{r}{E_n}{V_n}$ are the source and range maps, respectively. We will assume that $s^{-1}(v)\neq \varnothing$ for all $v\in V$ and that $r^{-1}(v)\neq \varnothing$ for all $v\in V\backslash V_0$. A Bratteli diagram
can be given a diagrammatic presentation with $V_n$ the vertices at level $n$ and $E_n$ the edges between $V_{n-1}$ and $V_n$. If $|V_{n-1}|=t_{n-1}$ and $|V_{n}|=t_{n}$ then the edge set $E_n$ is described by
a $t_n\times t_{n-1}$ incidence matrix $M_n=(m_{ij}^n)$, where $m_{ij}^n$ is the number of edges connecting $v_i^n\in V_n$ with $v_j^{n-1}\in V_{n-1}$ (see Figure \ref{fig:bD}). If the row sums are constant for every $M_n$, then we say that $(V,E)$ has the ERS-(Equal Row Sum) property. Let $k,l\in \mathbb{Z}^+$ with $k<l$ and let $E_{k+1}\circ E_k\circ \cdots
\circ E_l$ denote all the paths from $V_k$ to $V_l$.  Specifically, $E_{k+1}\circ E_k\circ \cdots \circ E_l = \set{(e_{k+1},\cdots,e_l)}{e_i\in E_i, i=k+1,\dots,l ; r(e_i)=s(e_{i+1}),i=k+1,\dots,l-1}$. We define $r\left((e_{k+1},\cdots,e_l)\right) = r(e_l)$ and $s\left((e_{k+1},\cdots,e_l)\right) = s(e_{k+1})$. Notice that the corresponding incidence matrix is the product $M_{l}M_{l-1}\cdots M_{k+1}$ of the individual incidence matrices.

\begin{definition}
\label{def:telescoping}
Given a Bratteli diagram $(V,E)$ and a sequence $0 = m_0 < m_1 < m_2 < \cdots$ in $\mathbb{Z^+}$, we define the \emph{telescoping} of $(V,E)$ to $\{m_n\}$ as $(V',E')$, where $V_n'=V_{m_n}$ and 
$E_n' = E_{m_{n-1}+1}\circ \cdots \circ E_{m_n}$, and the source and the range maps are as above.
\end{definition}

\begin{definition}
\label{def:finiterank}
 The Bratteli diagram $(V,E)$ is of \emph{finite rank} if $|V_n|\leq L <\infty$ for all $n$. By telescoping we may assume that $|V_n|=K$ for all $n=1,2,\dots$. We then say that $(V,E)$ is 
 of rank $K$, and write $\mathrm{rank}(V,E)=K$.
\end{definition}

\begin{definition}
\label{def:Simple}
We say that the Bratteli diagram $(V,E)$ is \emph{simple} if there exists a telescoping of $(V,E)$ such that the resulting Bratteli diagram $(V',E')$ has full connection between all consecutive levels, i.e.\ the entries of all the 
incidence matrices are non-zero.
\end{definition}

Given a Bratteli diagram $(V,E)$ we define the infinite path space associated to $(V,E)$, namely 
\[X_{(V,E)} = \set{(e_1, e_2, \dots)}{e_i\in E_i, r(e_i)=s(e_{i+1}); \quad \forall i \geq 1}.\]
Clearly $X_{(V,E)} \subseteq \prod_{n=1}^\infty E_n$, and we give $X_{(V,E)}$ the relative topology, $\prod_{n=1}^\infty E_n$ having the product topology. Loosely speaking this means that two paths in $X_{(V,E)}$ are close
if the initial parts of the two paths agree on a long initial stretch. Also, $X_{(V,E)}$ is a closed subset of $\prod_{n=1}^\infty E_n$, and so is compact. 

\medbreak
Let $p=(e_1, e_2, \dots, e_n)\in E_1\circ\cdots\circ E_n$ be a finite path starting at $v_0\in V_0$. We define the \emph{cylinder set} $U(p)=\set{(f_1,f_2,\dots)\in X_{(V,E)}}{f_i = e_i, i = 1,2,\dots,n}$. The collection of 
cylinder sets is a basis for the topology on $X_{(V,E)}$. The cylinder sets are clopen sets, and so $X_{(V,E)}$ is a compact, totally disconnected metric space. An admissable metric $d$ yielding 
the topology is $d(x,x')=\frac{1}{n}$ if $x=(e_1,e_2,\dots,e_{n-1},e_n,\dots)$, $y=(e_1,e_2,\dots,e_{n-1},e'_n,\dots)$, where $e_n\neq e'_n$.
If $(V,E)$ is simple then $X_{(V,E)}$ has no isolated points, and so $X_{(V,E)}$ is a Cantor set. (We will in the sequel disregard the trivial case where $|X_{(V,E)}|$ is finite.)

Let $P_n =E_1\circ\cdots\circ E_n$ be the set of finite paths of length $n$ (starting at the top vertex). We define the truncation map $\deff{\tau_n}{X_{(V,E)}}{P_n}$ by $\tau_n\left((e_1,e_2,\dots)\right) = (e_1,e_2,\dots,e_n)$.
If $m\geq n$ we have the obvious truncation map $\deff{\tau_{m,n}}{P_m}{P_n}$.

There is an obvious notion of isomorphism between Bratteli diagrams $(V,E)$ and $(V',E')$; namely, there exists a pair of bijections between $V$ and $V'$ preserving the gradings and intertwining the respective source and range
maps. Let $\sim$ denote the equivalence relation on Bratteli diagrams generated by isomorphism and telescoping. One can show that $(V,E)\sim(V',E')$ iff there exists a Bratteli diagram $(W,F)$ such that telescoping $(W,F)$ to odd
levels $0<1<3<\cdots$ yields a diagram isomorphic to some telescoping of $(V,E)$, and telescoping $(W,F)$ to even levels $0<2<4<\cdots$ yields a diagram isomorphic to some telescoping of $(V',E')$.

An \emph{ordered Bratteli diagram} $(V,E,\geq)$ is a Bratteli diagram $(V,E)$ together with a partial order $\geq$ in $E$ so that edges $e,e'\in E$ are comparable if and only if $r(e)=r(e')$. In other words, 
we have a linear order on each set $r^{-1}(v),\, v\in V\backslash V_0$. Assume $\abs{r^{-1}(v)}=m$ and the edge $e\in r^{-1}(v)$ has order $k$, where $1\leq k \leq m$. Then we will say that $e$ has 
\emph{ordinal} $k$, and we will write $\mathrm{ordinal}(e)=k$. We let $E_\text{min}$ and $E_\text{max}$, respectively, denote the minimal and maximal edges of the partially ordered set $E$.

Note that if $(V,E,\geq)$ is an ordered Bratteli diagram and $k<l$ in $\mathbb{Z}^+$, then the set $E_{k+1}\circ E_{k+2}\circ \cdots \circ E_l$ of paths from $V_k$ to $V_l$ with the same range can be given an induced
(lexicographic) order as follows: 
\[(e_{k+1}\circ e_{k+2}\circ \cdots \circ e_l)>(f_{k+1}\circ f_{k+2}\circ \cdots \circ f_l) \]
if for some $i$ with $k+1\leq i\leq l$, $e_j=f_j$ for $i<j\leq l$ and $e_i>f_i$. If $(V',E')$ is a telescoping of $(V,E)$ then, with this induced order from $(V,E,\geq)$, we get again an ordered Bratteli diagram $(V',E',\geq)$.

\begin{definition}
\label{def:properly}
 We say that the ordered Bratteli diagram $(V,E,\geq)$, where $(V,E)$ is a simple Bratteli diagram, is \emph{properly ordered} if
there exists a unique min path $x_\text{min} = (e_1,e_2,\dots)$ and a unique max path $x_\text{max} = (f_1,f_2,\dots)$ in $X_{(V,E)}$. (That is, $e_i\in E_\text{min}$ and $f_i\in E_\text{max}$ for all $i=1,2,\dots$.)
\end{definition}

Let $(V,E)$ be a properly ordered Bratteli diagram, and let $X_{(V,E)}$ be the path space associated to $(V,E)$. Then $X_{(V,E)}$ is a Cantor set. Let $T_{(V,E)}$ be the \emph{lexicographic map} on $X_{(V,E)}$, i.e.\ if 
$x=(e_1,e_2,\dots)\in X_{(V,E)}$ and $x\neq x_\text{max}$ then $T_{(V,E)}x$ is the successor of $x$ in the lexicographic ordering. Specifically, let $k$ be the smallest natural number so that $e_k\notin E_\text{max}$. Let $f_k$ be the 
successor of $e_k$ (and so $r(e_k)=r(f_k)$). Let $(f_1,f_2,\dots,f_{k-1})$ be the unique least element in $E_1\circ E_2 \circ \cdots \circ E_{k-1}$ from $s(f_k)\in V_{k-1}$ to the top vertex $v_0\in V_0$. 
Then $T_{(V,E)}((e_1,e_2,\dots)) = (f_1,f_2,\dots,f_k,e_{k+1},e_{k+2},\dots)$. We define $T_{(V,E)}x_\text{max}=x_\text{min}$. Then it is easy to check that $T_{(V,E)}$ is a minimal homeomorphism on $X_{(V,E)}$. We note that
if $x\neq x_\text{max}$ then $x$ and $T_{(V,E)}x$ are \emph{cofinal}, i.e.\ the edges making up $x$ and $T_{(V,E)}x$, respectively, agree from a certain level on. We will call the Cantor minimal system $(X_{(V,E)},T_{(V,E)})$ a
\emph{Bratteli-Vershik system}. There is an obvious way to telescope a properly ordered Bratteli diagram, getting another properly ordered Bratteli diagram, such that the associated Bratteli-Vershik systems are conjugate (cf.\ Definition \ref{def:factor})-- the map
implementing the conjugacy is the obvious one. By telescoping we may assume without loss of generality that the properly ordered Bratteli diagram has the property that at each level all the minimal edges (respectively the maximal edges)
have the same source, cf.\ \cite[Proposition 2.8]{HPS}.

We use the term \emph{dynamical system} to mean a compact metric space $X$ together with a homeomorphism $\deff{T}{X}{X}$, and we will denote this by $(X,T)$. We say $(X,T)$ is \emph{minimal} if
all $T$-orbits are dense. (Equivalently $T(A)=A$ for some closed $A\subseteq X$ implies that $A=X$ or $A=\varnothing$.) If $X$ is a Cantor set and $T$ is minimal, then we say that $(X,T)$ is a 
Cantor minimal system.

\begin{theorem}[\cite{HPS}]
\label{th:B-Vmod}
Let $(X,T)$ be a Cantor minimal system. Then there exists a properly ordered Bratteli diagram $(V,E,\geq)$ such that the associated Bratteli-Vershik system $(X_{(V,E)},T_{(V,E)})$ is conjugate to $(X,T)$. 
\end{theorem}

\begin{remark}
The simplest Bratteli-Vershik model $(V,E,\geq)$ for the odometer (see below) $(G_\mathfrak{a},T)$ associated to $\mathfrak{a}=(a_i)_{i\in\mathbb{N}}$ is obtained by letting $V_n=1$ for all $n$, and the number of edges between $V_{n-1}$ and
$V_n$ be $a_n$.
\end{remark}

Let $(G_\mathfrak{a},\rho_{\hat{\mathds{1}}})$ denote the \emph{odometer} (also called \emph{adding machine}) associated to the $\mathfrak{a}$-adic group
\[G_\mathfrak{a} = \prod_{i=1}^\infty \left\{0,1,\dots \frac{p_i}{p_{i-1}}-1\right\},\]
where $\mathfrak{a}= \left\{\frac{p_i}{p_{i-1}}\right\}_{i\in\mathbb{N}}$ (we set $p_0=1$) and where $\rho_{\hat{\mathds{1}}}(x) = x+\hat{1}$, where $\hat{1} = (1,0,0,\dots)$.
We note that $G_\mathfrak{a}$ is naturally isomorphic to the inverse limit group  
\[\quotient{\mathbb{Z}}{p_1\mathbb{Z}} \stackrel{\phi_1}{\longleftarrow}\quotient{\mathbb{Z}}{p_2\mathbb{Z}} \stackrel{\phi_2}{\longleftarrow} \quotient{\mathbb{Z}}{p_3\mathbb{Z}} \stackrel{\phi_3}{\longleftarrow} \cdots\]
where $\phi_i(n)$ is the residue of $n$ modulo $p_i$. It is a fact that the family consisting of compact groups $G$ that are both monothetic (i.e.\ contains a dense copy
of $\mathbb{Z}$ , which of course implies that $G$ is abelian) and Cantor (as a topological space), coincides with the family of 
$\mathfrak{a}$-adic groups. It is also noteworthy that all minimal rotations (in particular rotations by $\hat{1}$) on such groups are conjugate. This is a consequence of the fact that the 
dual group of an $\mathfrak{a}$-adic group is a torsion group.
If $\mathfrak{a} = \{p\}_{i\in\mathbb{N}}$, where $p$ is a prime, then $G_\mathfrak{a}$ is the $p$-adic integers. (We refer to \cite[Vol 1]{HR} for background information on $\mathfrak{a}$-adic groups.)

\begin{remark}
 It is well known, and easy to prove, that the Cantor minimal system $(X,T)$ is conjugate (cf.\ Definition \ref{def:factor}) to an odometer if and only if it is the inverse limit of a sequence
 of periodic systems.
\end{remark}

\begin{definition}
\label{def:expansive}
 $(X,T)$ is \emph{expansive} if there exists $\delta>0$ such that if $x\neq y$ then $\mathrm{sup}_n \mathrm{d}(T^nx, T^n y)>\delta$, where $\mathrm{d}$ 
is a metric that gives the topology of $X$. (Expansiveness is independent of the metric as long as the metric gives the topology of $X$.)
\end{definition}

Let $\Lambda=\{a_1,a_2,\dots,a_n\},\, n\geq 2$, be a finite alphabet and let $Z=\Lambda^\mathbb{Z}$ be the set of all bi-infinite sequences of symbols from
$\Lambda$ with $Z$ given the product topology -- thus $Z$ is a Cantor set. Let $\deff{S}{Z}{Z}$ denote the shift map, $\deff{S}{(x_n)}{(x_{n+1})}$.
If $X$ is a closed subset of $Z$ such that $S(X)=X$, we say that $(X,S)$ is a \emph{subshift}, where we denote the restriction of $S$ to $X$ again by $S$.
Subshifts are easily seen to be expansive. We state the following well-known fact as a proposition. (Cf.\ \cite[Theorem 5.24]{walters}.)

\begin{proposition}
\label{prop:subshift}
Let $(X,T)$ be a Cantor minimal system. Then $(X,T)$ is conjugate to a minimal subshift on a finite alphabet if and only if $(X,T)$ is expansive.
\end{proposition}

\begin{definition}
\label{def:factor}
We say that a dynamical system $(Y,S)$ is a \emph{factor} of $(X,T)$ and that $(X,T)$ is an \emph{extension} of $(Y,S)$ if there exists a continuous surjection $\deff{\pi}{X}{Y}$ which 
satisfies $S(\pi(x))=\pi(Tx),\,\forall x\in X$. We call $\pi$ a \emph{factor map}.
If $\pi$ is a bijection then we say that $(X,T)$ and $(Y,S)$ are \emph{conjugate}, and we write $(X,T)\cong (Y,S)$. 
\end{definition}

Let $(V,E,\geq)$ be a properly ordered Bratteli diagram, and let $(X_{(V,E)},T_{(V,E)})$ be the associated Bratteli-Vershik system. For each $k\in \mathbb{N}$ let $P_k$ as above denote the paths from $V_0$ to $V_k$, i.e.\ the paths from $v_0\in V_0$ 
to some $v\in V_k$. For $x\in X_{(V,E)}$ we associate the bi-infinite sequence $\pi_k(x) = \left(\tau_k(T_{(V,E)}^n x)\right)_{n=-\infty}^\infty\in P_k^\mathbb{Z}$ over the finite alphabet $P_k$, where $\deff{\tau_k}{X_{(V,E)}}{P_k}$ 
is the truncation map. Let $S_k$ denote the shift map on $P_k^\mathbb{Z}$. Then the following diagram commutes
\begin{center}
 \includegraphics{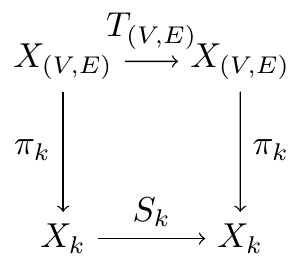}
\end{center}
where $X_k = \pi_k(X_{(V,E)})$. One observes that $\pi_k$ is a continuous map, and so $X_k$ is a compact shift-invariant subset of $P_k^\mathbb{Z}$. So $(X_k,S_k)$ is a factor of $(X_{(V,E)}, T_{(V,E)})$. For $k>l$ there is an 
obvious factor map $\deff{\pi_{k,l}}{X_k}{X_l}$, and one can show that $(X_{(V,E)}, T_{(V,E)})$ is the inverse limit of the system $\{(X_k,S_k)\}_{k\in\mathbb{N}}$.
We write $(X_{(V,E)},T_{(V,E)})=\varprojlim (X_k,S_k)$. All the systems $(X_k,S_k)$ are clearly expansive. One has the following result which will be important for us.

\begin{proposition}
Assume $(X_{(V,E)}, T_{(V,E)})$ is expansive. Then there exists $k_0\in \mathbb{N}$ such that for all $k\geq k_0$, $(X_{(V,E)}, T_{(V,E)})$ is conjugate to $(X_k,S_k)$ by the map $\deff{\pi_k}{X_{(V,E)}}{X_k}$.
\end{proposition}

\proc{Proof.}
Since the $\pi_k$'s are factor maps, all we need to show is that there exists $k_0$ such that $\pi_k$ is injective for all $k\geq k_0$. Recall that $(X_{(V,E)}, T_{(V,E)})$ being expansive means that there exists $\delta>0$ such that given 
$x\neq y$ there exists $n_0\in \mathbb{Z}$ such that $d(T_{(V,E)}^{n_0}x,T_{(V,E)}^{n_0}y)>\delta$, where $d$ is some metric on $X_{(V,E)}$ compatible with the topology. Choose $k_0$ such that $d(x,y)<\delta$ if $x$ and $y$ agree (at least) on the $k_0$ first edges. Now assume that 
$\pi_k(x)=\pi_k(y)$ for some $k\geq k_0$. By the definition of $\pi_k$ this means that, for all $n\in \mathbb{Z}$, $\tau_k(T_{(V,E)}^n x) = \tau_k(T_{(V,E)}^n y)$, and so $d(T_{(V,E)}^n x,T_{(V,E)}^n y)<\delta$ for all $n\in\mathbb{Z}$ because of our choice of $k_0$. 
This contradicts that $d(T_{(V,E)}^{n_0}x,T_{(V,E)}^{n_0}y)>\delta$. Hence $\pi_k$ is injective for all $k\geq k_0$, proving the proposition.
$\square$
\medbreak

We draw the following conclusions from the above: Let $(X_{(V,E)},T_{(V,E)})$ be the Bratteli-Vershik system associated to the properly ordered Bratteli diagram $(V,E,\geq)$. Then $(X_{(V,E)},T_{(V,E)})$
is not expansive if and only if $\deff{\pi_k}{X_{(V,E)}}{X_k\, (=\pi_k(X_{(V,E)}))}$ is not injective for $k=1,2,3,\dots$.

\section{Key definitions and basic properties.}
\label{sec:defi}

Set $X=X_{(V,E)}$, $T=T_{(V,E)}$, where $(X_{(V,E)},T_{(V,E)})$ is the Bratteli-Vershik system associated to the properly ordered Bratteli diagram $(V,E,\geq)$. (We will use the notation introduced
in Section \ref{sec:bratteli} as well as the one in \cite{DM}, and we adopt the terminology of \cite{DM}.)

Consider a pair $(x,x')$ of distinct points in $X$ such that $\pi_i(x)=\pi_i(x')$ for some $i\geq 1$. We call such a pair \emph{$i$-compatible}. Observe that $(x,x')$ is then $k$-compatible if $k\leq i$. Since $x\neq x'$, there exists some $j> i$ such that 
$\pi_j(x)\neq \pi_j(x')$. We say that the pair is \emph{$j$-separated}. The largest index $i_0$ for which the pair $(x,x')$ is $i_0$-compatible (and hence it is $(i_0+1)$-separated) will be
called the \emph{depth of compatibility} (\emph{depth} for short) of this pair. In particular, equal elements have depth $\infty$. Let $(x,x')$ be $i$-compatible and $j$-separated for some $j>i$.
By telescoping between levels $i$ and $j$ we obtain that $(x,x')$ is of depth $i$,  which is easily seen.
\newline

We make some observations:
\begin{enumerate}[(i)]
 \item \label{Pt:1} If $(x,x')$ is $i$-compatible and $j$-separated, then $(T^m x,T^m x')$ is $i$-compatible and $j$-separated for all $m\in\mathbb{Z}$. [This follows since $\pi_k(T^my)=S_k^m\pi_k(y)$ for 
 all $y\in X$, $k=1,2,3,\dots$.]
 \item \label{Pt:2} If $(x,x')$ is of depth $i$, then $(T^mx,T^mx')$ is of depth $i$ for all $m\in\mathbb{Z}$. [This is an immediate consequence of \eqref{Pt:1}.]
 \item \label{Pt:3} If $(x,x')$ is a pair of depth $i$ and $(x,x'')$ is a pair of depth $j>i$, then $(x',x'')$ is a pair of depth $i$ (and hence not equal). [Clearly the pair $(x',x'')$ is $i$-compatible. There
 exists $m\in\mathbb{Z}$ such that $\tau_{i+1}(T^mx)\neq \tau_{i+1}(T^mx')$. Since $\tau_{i+1}(T^mx)= \tau_{i+1}(T^mx'')$, the assertion follows.]
\end{enumerate}

An $i$-compatible and $j$-separated ($j>i$) pair $(x,x')$ is said to have a \emph{common $j$-cut} if for some $m\in\mathbb{Z}$, $\tau_j(T^mx)$ and $\tau_j(T^mx')$ are minimal paths, i.e.\ consisting of only minimal
edges, between level $j$ and level 0 (i.e.\ the top vertex). Note that if a pair has a common $j$-cut it also has a common $j'$-cut for every $i< j'\leq j$. It is obvious from the definitions that 
if $(x,x')$ has a common j-cut, then $(T^lx,T^lx')$ also has a common $j$-cut for any $l\in\mathbb{Z}$. Observe also that if the pair $(x,x')$ has no common $j$-cut the pair must be $j$-separated.

\begin{figure}[h]
\centering
 \includegraphics{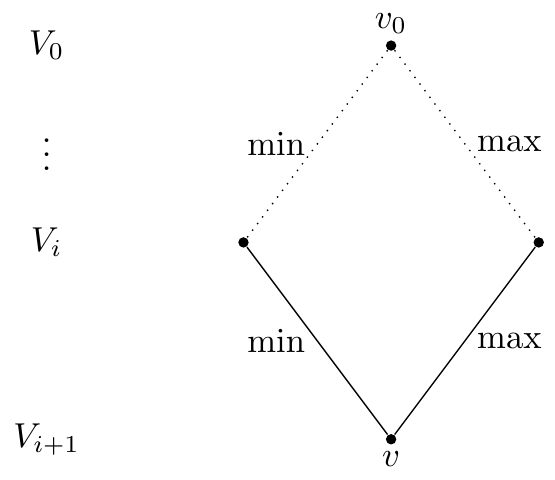}
\caption{}
\label{fig:separated}
\end{figure}
We make one important observation: Let $(x,x')$ be of depth $i$, and assume $(x,x')$ has a common $(i+1)$-cut. Then for some $m\in\mathbb{Z}$, the pair $(T^mx,T^mx')$ is of depth $i$ such that 
$\tau_{i+1}(T^mx)$ and $\tau_{i+1}(T^mx')$ are minimal paths, \emph{and} $r(\tau_{i+1}(T^mx))\neq r(\tau_{i+1}(T^mx'))$. [In fact, by assumption there exists $k\in\mathbb{Z}$ such that $\tau_{i+1}(T^kx)$ and
$\tau_{i+1}(T^kx')$ are minimal paths. If $v= r(\tau_{i+1}(T^kx))= r(\tau_{i+1}(T^kx'))$ then $l$ iterates of $T$, say, applied to $T^kx$ and $T^kx'$ respectively, will ``sweep over'' all the paths between
$v_0\in V_0$ and $v\in V_{i+1}$, eventually reaching the max path, see Figure \ref{fig:separated}. Applying $T$ one more time to $T^{k+l}x$ and $T^{k+l}x'$, respectively, will result in $\tau_{i+1}(T^px)$ and 
$\tau_{i+1}(T^px')$ are minimal paths. (Here $p=k+l+1$.) If $r(\tau_{i+1}(T^px))\neq r(\tau_{i+1}(T^px'))$ we are done, setting $m=p$. If $r(\tau_{i+1}(T^px))= r(\tau_{i+1}(T^px'))$, we do the 
same procedure as above. If we get to a stage where the ranges are distinct we are done. If this does not happen, we play the same game on $T^kx$ and $T^kx'$, but now with iterates of $T^{-1}$ instead
of $T$. This must lead to a situation where the ranges are distinct, otherwise $\pi_{i+1}(x)=\pi_{i+1}(x')$, contradicting that $(x,x')$ is $(i+1)$-separated.] 

\section{Proof of THEOREM.}
We assume that $(X_{(V,E)},T_{(V,E)})$ it not expansive and so for all $i\geq 1$, $\deff{\pi_i}{X}{X_i}$ is not injective. This is easily seen to have as a consequence that for infinitely many levels $i$
there exist pairs of points $(x_i,x_i')$ of depth $i$. If we telescope between these levels we may assume that for every $i\geq 1$ there exists a pair $(x_i,x_i')$ of depth $i$. We will show that $(X_{(V,E)},T_{(V,E)})$ is an odometer, which will 
complete the proof. First we set the stage in the sense that we may assume that $(V,E,\geq)$ has the following properties:
\begin{enumerate}[(i)]
 \item \label{pt:en} We may assume that $\mathrm{rank}(V,E)=K$ (cf.\ Definition \ref{def:finiterank}) is the smallest possible such that the Bratteli-Vershik system associated to $(V,E,\geq)$ is (conjugate to) the given one.
 (If $K=1$ we have an odometer, so there is nothing more to prove.)
 \item \label{pt:to} By  telescoping we may assume that between consecutive levels there is full connection (cf.\ Definition \ref{def:Simple}) and, furthermore, that at each level all the minimal
 edges (respectively the maximal edges) have the same source. (This is not an essential assumption, but it makes it easier to visualize the Vershik map.) 
\end{enumerate}
Note that the property \eqref{pt:en} is not affected by the operations performed in \eqref{pt:to}.
\newline

As before we let $X=X_{(V,E)}$, $T=T_{(V,E)}$. There are two scenarios, mutually exclusive, cf.\ \cite{DM}.
\begin{enumerate}[(1)]
 \item \label{pt:cut} There exists $i_0$ such that for all $i\geq i_0$ and every $j>i$ there exists a pair $(x,x')$ of depth $i$ with a common $j$-cut.
 \item \label{pt:nocut} For infinitely many $i$, any pair $(x,x')$ of depth $i$ has no common $j$-cuts for sufficiently large $j>i$. (Note that $j$ depends upon $(x,x')$!)
\end{enumerate}
The proof is different for case \eqref{pt:cut} and case \eqref{pt:nocut}. 
\newline

We consider  case \eqref{pt:cut}:
\newline

\begin{figure}[h]
\centering
 \includegraphics[scale=1]{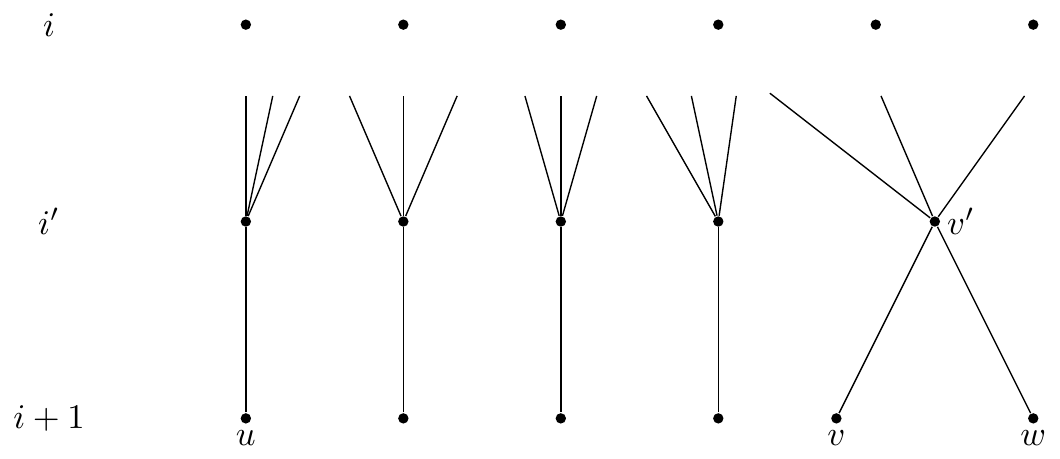}
\caption{Rank $K=6$.}
\label{fig:reduce}
\end{figure}

\begin{figure}[h]
\centering
 \includegraphics[scale=1]{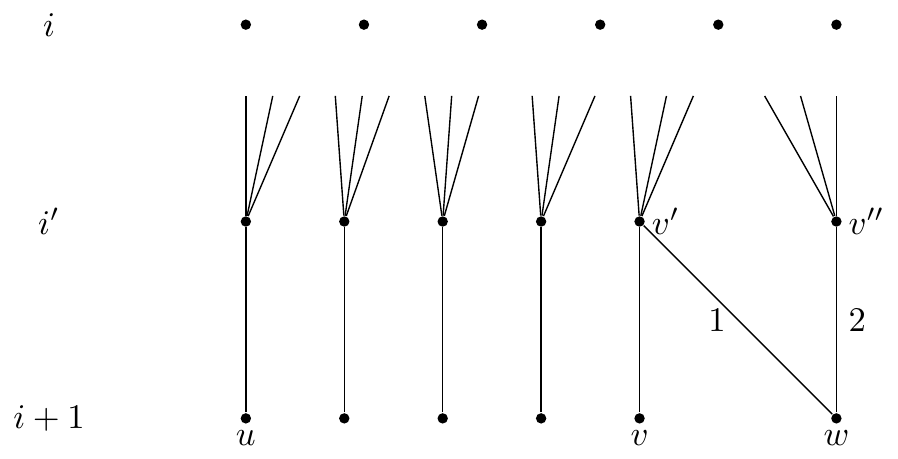}
\caption{Rank $K=6$.}
\label{fig:duplicate}
\end{figure}

By telescoping we may assume that for each $i\geq i_0$ there exists a pair $(x,x')$, $x,x'\in X$, of depth $i$. The idea is to find another properly ordered Bratteli diagram $(V',E',\geq)$ with 
$\mathrm{rank}(V',E')<K$ (assuming $K>1$), such that $(X_{(V',E')},T_{(V',E')})\cong (X,T)$. This contradiction will finish the proof in this case. Now choose any $i\geq i_0$. By the observation we made at the end of 
Section \ref{sec:defi} we may assume that there exists a pair $(x,x')$ of depth $i$ such that $\tau_{i+1}(x)$ and $\tau_{i+1}(x')$ consist of minimal edges, and that $v=r(\tau_{i+1}(x))\neq r(\tau_{i+1}(x'))=w$.
If $|r^{-1}(v)|=|r^{-1}(w)|$ we may insert a new level (we name it $i'$) between levels $i$ and $i+1$ with ordering of the edges as shown in Figure \ref{fig:reduce}. 
(The ordering at the vertex $v'$ is the same as the ordering at $v$ and $w$, the two latter being the same since $(x,x')$ is of depth $i$.) The order of the edges ranging at vertices $u\in V_{i+1}-\{v,w\}$ is replicated at level $i'$.
We notice that if we telescope between levels $i$ and $i+1$ we get the original ordering. So the insertion of level $i'$ does not change the Bratteli-Vershik map. Now we have obtained a level $i'$
with $K-1$ number of vertices. If $|r^{-1}(v)|<|r^{-1}(w)|$, say, we insert a new level $i'$ between levels $i$ and $i+1$ as shown in Figure \ref{fig:duplicate}. The $|r^{-1}(v)|$ first edges ranging 
at $v$ and $w$ are ordered at $v'$ as they are at $v$ and $w$ while the $|r^{-1}(w)|-|r^{-1}(v)|$ remaining edges ranging at $w$ are ordered at $v''$ as they are at $w$. As before vertices $u\in V_{i+1}-\{v,w\}$ 
are just replicated at level $i'$. We observe that the number of vertices at level $i'$ is the same as at level $i+1$, namely $K$. Now we claim that $(x,x')$ separates at level $i'$, and so $(x,x')$ 
has depth $i$ in the new diagram as well. In fact, by applying $L+1$ iterates of $T$ to $x$ and $x'$, respectively, we see that they separate at level $i'$.  Here $L$ is the number of paths from the top 
vertex ranging at $v$. We observe that the number of edges between levels $i$ and $i'$ is strictly smaller than the number of edges between levels $i$ and $i+1$. Now we repeat the same construction between levels $i$ and $i'$. Since we decrease the number of 
edges each time, we must eventually arrive at the first case, where the number of edges ranging at $v$ and $w$ are the same. Doing the construction we did in the first case will then yield a level which has $K-1$ vertices. 

After we have done this, we do the same construction between levels $(i+1)$ and $(i+2)$, etc. If we now telescope to the new levels with $K-1$ vertices we wind up with a properly ordered Bratteli
diagram of rank $K-1$ which yields a Bratteli-Vershik system conjugate to the original. From this we conclude that $K$ can not be larger than one and the proof is completed for case \eqref{pt:cut}.
\newline

We now look at case \eqref{pt:nocut}. By telescoping to appropriate levels we may assume that we have the following scenario:
\newline

For each $i\geq 1$ there exists a pair of depth $i$ that has no common $(i+1)$-cuts, and hence no common $j$-cuts for any $j>i$. Now fix any $i_0\geq1$. We shall prove that $(X_{i_0},S_{i_0})$ is periodic (and hence
finite). This will imply that $(X,T)$ is an odometer since $(X,T)=\varprojlim (X_k,S_k)$, hence finishing the proof.

Under the assumption that the above scenario holds we can prove the following lemma (whose proof we postpone to the end).

\begin{lemma*}
 For any positive integer $L$ there exist $L$ distinct elements $y_1,y_2, \dots,y_{L}$ in $X$ which are pairwise $i_0$-compatible and pairwise have no common $j$-cuts for some 
 $j\geq i_0$. In particular, they are pairwise $j$-separated. (Observe that for all $k\in\mathbb{Z}$ the elements $T^ky_1,T^ky_2,\dots,T^ky_L$ have the same properties as
 $y_1, y_2, \dots, y_L$.)
\end{lemma*}

By telescoping between level $i_0$ and level $j$ we may assume that the elements $y_1, y_2,\dots,y_L$ in the Sublemma are pairwise of depth $i_0$ and have no common $(i_0+1)$-cuts. Choose $L$ in the 
Sublemma to be \[L=(K-1)2^K+2.\]

\begin{figure}[h]
\centering
 \includegraphics[width=\textwidth]{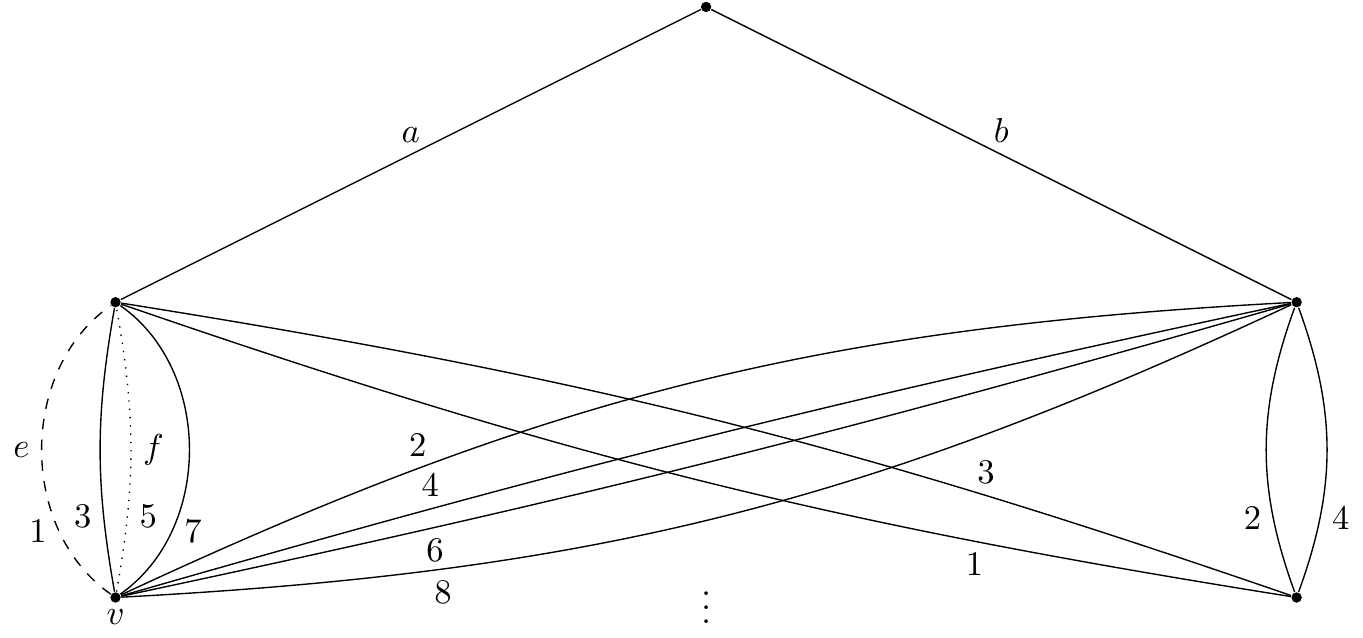}
\caption{A Bratteli diagram where $i_0=1$, $\abs{v}=8$, $l_v=4$, $y_1=(a,e,\dots)$, $y_2=(a,f,\dots)$ and $\pi_1(y_1)=\pi_1(y_2)$.}
\label{fig:periodic}
\end{figure}

Let us in the sequel denote $\tau_{i_0}$ by $\tau_1$ and $\tau_{i_0+1}$ by $\tau_2$. Let $v\in V_{i_0+1}$ and let $l_v$ be the smallest (positive) difference of the ordinal numbers of any pair $(\tau_2(T^py_i),\tau_2(T^py_j))$ with common range $v$, i.e.\ $r(\tau_2(T^py_i))=r(\tau_2(T^py_j))=v$. Here $i,j\in \{1,2,\dots,L\}$, $i\neq j$, and $p$ can be any integer. (In Figure \ref{fig:periodic} we have illustrated this by assuming that $l_v$ is obtained at $v$
by $y_1=(a,e,\dots)$, $y_2=(a,f,\dots)$. We see that $l_v=5-1=4$. Actually, Figure \ref{fig:periodic} illustrates another point (setting aside that $K=2$): In the general case, if we telescope between level 0 and level $i_0$, 
then we wind up with a scenario like the one in Figure \ref{fig:periodic} except that there are multiple edges instead of the single edge $a$ (respectively $b$).)

Assume $l_v$ is obtained at $v$ with the pair $(\tau_2(T^py_i),\tau_2(T^py_j))$. Since \\$\pi_{i_0}(T^py_i)=\pi_{i_0}(T^py_j)$ this has the following consequence: $\hat{v}(k) = \hat{v}(k+l_v)$ for 
$k\in [1,\abs{v}-l_v]$. Here $\abs{v}$ denotes the number of paths in $P_{i_0+1}$ ranging at $v$, and $\hat{v}(k)$ is the element in $P_{i_0}$ obtained by ``cutting off'' (or truncating) the path in $P_{i_0+1}$ 
ranging at $v$ with ordinal number $k$. (In Figure \ref{fig:periodic} we have $l_v=4$, and we get that $\hat{v}(1)=\hat{v}(1+4)=a$, $\hat{v}(2)=\hat{v}(2+4)=b$, $\hat{v}(3)=\hat{v}(3+4)=a$, 
$\hat{v}(4)=\hat{v}(4+4)=b$.)

Now let $\hat{y}$ denote the common image of $y_1, y_2,\dots,y_L$ under $\pi_{i_0}$, i.e.\ $\pi_{i_0}(y_1)= \pi_{i_0}(y_2)=\cdots=\pi_{i_0}(y_L)=\hat{y}\subseteq P_{i_0}^\mathbb{Z}$. 
Observe that by the definition of $\hat{y}$ we have that $\hat{y}(l)=\tau_1(T^ly_i)$ for $l\in\mathbb{Z}$ and any $i=1,2,\dots,L$. In particular, $\hat{y}(0)=\tau_1(y_i)$. We will say that the $l_v$-periodicity
law holds at the coordinate $n\in\mathbb{Z}$ of $\hat{y}$ if $\hat{y}(n)=\hat{y}(n+l_v)$. We make one important observation: If, say, $r(\tau_2(T^ny_i))=r(\tau_2(T^ny_j))=v$ for some $y_i\neq y_j$
and $\tau_2(T^ny_i)<\tau_2(T^ny_j)$ ($n\in\mathbb{Z}$), then the $l_v$ periodicity law holds at the coordinate $n$ of $\hat{y}$. In fact, if the ordinal number of $\tau_2(T^n y_i)$ is $k$, then
$k+l_v\leq$ (ordinal number of $\tau_2(T^n y_j)$) $ \leq \abs{v}$, and so $\hat{v}(k)=\hat{v}(k+l_v)$. By definition of $\hat{y}$ it follows that $\hat{v}(k)=\hat{y}(n)$. 
Now $\tau_2(T^{n+l_v}y_i)\leq \tau_2(T^{n}y_j)$, and so $\hat{y}(n+l_v)=\hat{v}(k+l_v)$, and hence the $l_v$-periodicity law holds at the coordinate $n$ of $\hat{y}$.

Let us order the vertices at level $i_0+1$ by $v_1,v_2,\dots,v_K$ such that $l_{v_1}\leq l_{v_2}\leq \cdots \leq l_{v_K}$. (If there exists some vertex $v$ at level $i_0+1$ such that no two 
$T^py_i, T^py_j$ ($p\in\mathbb{Z}, i\neq j$) range at $v$, then we just ignore that $v$. This will not cause any problem for the subsequent argument, so we may just as well assume that there exists 
no such $v$.) Assume that there exists some $k\in\mathbb{Z}$ such that 
\[\abs{\{T^ky_i\left|r(\tau_2(T^ky_i))=v_1,\, i=1,2,\dots,L\right.\}}\geq K+1.\]
(Here $\abs{A}$ denotes the cardinality of the set $A$.) By renaming $T^ky_i$ as $y_i$, $i=1,2,\dots,L$ (cf.\ Sublemma), we may assume 
\[\abs{\{y_i\left|r(\tau_2 (y_i))=v_1,\, i=1,2,\dots,L\right.\}}\geq K+1.\]
Let $I$ be the largest interval of integers (obviously containing 0) such that the $l_{v_1}$-periodicity law holds. Specifically, if $i\in I$, then $\hat{y}(i)=\hat{y}(i+l_{v_1})$. If $I$ is 
infinite at the right end, 
then $\hat{y}(i)=\hat{y}(i+l_{v_1})$ for all $i\geq 0$. Shifting $\hat{y}$ to the left and using minimality of $(X_{i_0},S_{i_0})$, we get that $\hat{y}$ is periodic and so $(X_{i_0},S_{i_0})$ is 
periodic, thus finishing the proof. 
If $I$ has a right end, let $m\in\mathbb{Z}_+$ be the first integer to the right of $I$. At least two of the elements in $\{T^my_i\,|\, r(\tau_2(y_i))=v_1, i=1,2,\dots,L\}$, say, $T^my_i$ and 
$T^my_j$ ($i\neq j$) are such that $r(\tau_2(T^my_i))=r(\tau_2(T^my_j))=v_k$ for some $k\geq 1$.
We have $\hat{y}(m)=\tau_1(T^my_j)=\tau_1(T^my_i)$. If $v_k=v_1$ then the $l_{v_1}$-periodicity law holds at $m$ by the observation we made above, contradicting our assumption. So $k>1$. Let the 
ordinal numbers of $\tau_2(T^my_i)$ and $\tau_2(T^my_j)$ be $s$ and $t$, respectively, and assume $s<t$. Now $l_{v_1}\leq l_{v_k}\leq t-s$, and so the ordinal number $s+l_{v_1}$ exists for paths in $P_{i_0+1}$ 
ranging at $v_k$.

Applying $T^{-(t-s)}$ to $T^my_j$ results in the following:
\[\tau_2(T^my_i)=\tau_2(T^{-(t-s)}(T^my_j))=\tau_2(T^{m-(t-s)}y_j).\]
In particular, the ordinal numbers of $\tau_2(T^{m-(t-s)}y_j)$ and $\tau_2(T^{m}y_i)$ are the same, both equal to $s$. We also get 
\[\hat{y}(m)=\tau_1(T^my_i)=\tau_1(T^{m-(t-s)}y_j)=\hat{y}(m-(t-s)).\]
Assume we can prove that $m-(t-s)>0$. Then $m-(t-s)\in I$ and so the $l_{v_1}$ periodicity law holds at $m-(t-s)$ i.e.\ $\hat{y}(m-(t-s)+l_{v_1})=\hat{y}(m-(t-s))$. If we apply $T^{l_{v_1}}$ to both $T^my_i$ and $T^{m-(t-s)}y_j$, respectively,
we get \[\tau_2(T^{m+l_{v_1}}y_i)=\tau_2(T^{m-(t-s)+l_{v_1}}y_j)\]
(both having ordinal number $s+l_{v_1}$), and so
\[\hat{y}(m+l_{v_1})=\hat{y}(m-(t-s)+l_{v_1})=\hat{y}(m-(t-s))=\hat{y}(m).\]
So the $l_{v_1}$-periodicity law holds at $m$ which contradicts our assumption that $I$ has a (finite) right end, thus finishing the proof. It remains to prove that $m-(t-s)>0$. Assume by contradiction that $m-(t-s)\leq 0$. There exists
an $l$, $0\leq l \leq t-s$ ($\leq t$) such that $m-l=0$. We get
\[\tau_2(y_j)=\tau_2(T^{m-l}y_j)=\tau_2(T^{-l}(T^my_j))\]
which is impossible since $r(\tau_2(y_j))=v_1$ and $r(\tau_2(T^{-l}(T^my_j)))=v_k$.

To recap, if there exists some $k\in\mathbb{Z}$ such that 
\begin{equation}
 \label{eq:ranges}
 \abs{\{T^ky_i\,|\, r(\tau_2(T^ky_i))=v_1,\, i=1,2,\dots,L\}}\geq K+1,\tag{*}
\end{equation}
then we can prove that $(X_{i_0},S_{i_0})$ is periodic. So assume that this is not the case. In other words, for all $k\in\mathbb{Z}$ we have 
\[ \abs{\{T^ky_i\,|\, r(\tau_2(T^ky_i))=v_1,\, i=1,2,\dots,L\}}\leq K.\]
Assume now that there exists some $k\in\mathbb{Z}$ such that 
\begin{equation}
 \label{eq:ranges2}
 \abs{\{T^ky_i\,|\, r(\tau_2(T^ky_i))=v_2,\, i=1,2,\dots,L\}}\geq 2K.\tag{**}
\end{equation}
Now we argue exactly as above letting $I$ and $m$ be as above. There will then exist at least two elements in $\{T^my_i\,|\, r(\tau_2(y_i))=v_2,\, i=1,2,\dots,L\}$, say $T^my_i$ and $T^my_j$ ($i\neq j$) such that $r(\tau_2(T^my_i))=r(\tau_2(T^my_j))=v_k$, where $k\geq 2$. By exactly 
the same argument as above, we get that $(X_{i_0},S_{i_0})$ is periodic. If both \eqref{eq:ranges} and \eqref{eq:ranges2} do not occur, we assume there exists $k\in\mathbb{Z}$ such that
\begin{equation}
 \label{eq:ranges3}
 \abs{\{T^ky_i\,|\, r(\tau_2(T^ky_i))=v_3,\, i=1,2,\dots,L\}}\geq 4K-2.\tag{***}
\end{equation}
We repeat the same argument as above, again getting that $(X_{i_0},S_{i_0})$ is periodic. We continue this process, and it must eventually stop. The ``worst'' case scenario is that for all $k\in\mathbb{Z}$ the following simultaneously holds:
\begin{equation}
\label{eq:sum}
 \left.\begin{aligned}
\abs{\{T^ky_i\,|\, r(\tau_2(T^ky_i))=v_1,\, i=1,2,\dots,L\}}&\leq K\\
\abs{\{T^ky_i\,|\, r(\tau_2(T^ky_i))=v_2,\, i=1,2,\dots,L\}}&\leq 2K-1\\
\abs{\{T^ky_i\,|\, r(\tau_2(T^ky_i))=v_3,\, i=1,2,\dots,L\}}&\leq 4K-3\\
\vdots \hspace{1.5cm}\vdots \hspace{4cm} & \hspace{1cm} \vdots \\
\abs{\{T^ky_i\,|\, r(\tau_2(T^ky_i))=v_l,\, i=1,2,\dots,L\}}&\leq 2^{l-1}K-(2^{l-1}-1)\\
\vdots \hspace{1.5cm}\vdots \hspace{4cm} & \hspace{1cm} \vdots \\
\abs{\{T^ky_i\,|\, r(\tau_2(T^ky_i))=v_K,\, i=1,2,\dots,L\}}&\leq 2^{K-1}K-(2^{K-1}-1)
       \end{aligned}
 \right\} \tag{****}
\end{equation}

Adding up the right hand side we get $(K-1)2^K+1$. However, \eqref{eq:sum} contradicts that $L=(K-1)2^K+2$. Hence one of the scenarios that lead to a proof that $(X_{i_0},S_{i_0})$ is periodic must occur. This
finishes the proof of THEOREM. \qed

\begin{proof}[Proof of Sublemma]
 For each $i\in[i_0,i_0+L-1]$ let $(x_i,x_i')$ be a pair of depth $i$ which do not have a common $(i+1)$-cut. Let $\{n_k\}_k$ be a subsequence of natural numbers such that $T^{n_k}x_i\longrightarrow y_0$ as 
 $k\longrightarrow \infty$, where $y_0$ is the unique minimal path $x_\mathrm{min}$ in $X=X_{(V,E)}$. (Because of minimality of $(X,T)$ such a subsequence exists.) By compactness of $X$ there exists a subsequence 
 of $\{n_k\}_k$, which we again will denote by $\{n_k\}_k$, such that $T^{n_k}x_i'\longrightarrow y_i$ for some $y_i\in X$. By continuity we get that $\pi_i(y_0)=\pi_i(y_i)$ since $\pi_i(T^{n_k}x_i)=\pi_i(T^{n_k}x_i')$
 for all $k$. So $(y_0,y_1)$ is $i$-compatible. We claim that $(y_0,y_i)$ is $(i+1)$-separated, and hence $(y_0,y_i)$ is of depth $i$. In fact, there exists $k_0$ such that for all $k\geq k_0$, $\tau_{i+1}(T^{n_k}x_i)=\tau_{i+1}(y_0)$ and 
 $\tau_{i+1}(T^{n_k}x_i')=\tau_{i+1}(y_i)$. Since $(x_i,x_i')$ do not have a common $(i+1)$-cut we conclude that $\tau_{i+1}(T^{n_k}x_i)\neq \tau_{i+1}(T^{n_k}x_i')$. Hence $\tau_{i+1}(y_0)\neq \tau_{i+1}(y_i)$,
 and so $(y_0,y_i)$ is $(i+1)$-separated, hence of depth $i$. By \eqref{Pt:3} in Section \ref{sec:defi} we get that if $i<i'$, then $(y_i,y_{i'})$ is a pair of depth $i'$ and hence, in particular, $i_0$-compatible.
 (In particular, the points $y_{i_0}, y_{i_0+1}, \dots, y_{i_0+L-1}$ are distinct.) By assumption \eqref{pt:nocut} there exists $j(i,i')>i'$ such that $(y_i,y_{i'})$ have no common $j(i,i')$-cut. Let 
 $j=\mathrm{max}\{j(i,i')|i\neq i'\}$. Then $y_1, y_2, \dots , y_L$ are distinct points in $X$ which are pairwise $i_0$-compatible and pairwise have no common $j$-cut. (Here we rename the indices by letting 
 $i_0\rightarrow 1, i_0+1\rightarrow 2, \dots, i_0+L-1\rightarrow L$.) This finishes the proof of the Sublemma.
\end{proof}

The following corollary gives a positive answer to QUESTION 1 raised in \cite{DM} about finding a smaller $L$ than the one given in the so-called "Infection Lemma" in \cite{DM}, namely $L=K^{K+1}+1$. In fact,
we conjecture that the $L$ we have found is optimal.

\begin{corollary}
\label{cor:bound}
If there exists at least \[L=(K-1)2^K+2\]
points $y_k$ ($k\in[1,L]$) that are $i$-compatible and have no common $j$-cut for some $j>i$ (and hence are $j$-separated), then $(X_{i},S_{i})$ is periodic.
\end{corollary}

\begin{remark}
We find some of the assertions at the beginning of "Proof in case (2)" of THEOREM 1 (the same as our THEOREM) in \cite[pp. 744-745]{DM} somewhat confusing and lacking explanations. For example, it is stated that under assumption (2) 
the following holds: for each $i\geq 1$ every (!) pair of depth $i$ has no common $(i+1)$-cuts. We do not see why this should be true. However, it does follow from assumption (2) that there exists (!) a pair with the desired properties,
and that is sufficient for the proof to work. We also find the subsequent argument for the existence of appropriately many $i_0$ compatible and $j$-separated elements with no common $j$-cuts somewhat lacking and needing some further 
explanation.

While we think that there is too much "hand-waving" in the proof of THEOREM 1 in \cite{DM}, making it difficult to follow, we want no negative impression to attach to the \cite{DM} paper. In fact, Downarowicz and Maass had the insight 
to realize that such a remarkable result as THEOREM 1 holds, and also the ingenuity of finding a proof, the basic idea of which we use in our new proof.
\end{remark}

\addcontentsline{toc}{section}{Bibliography}
\bibliographystyle{amsalpha}

\end{document}